\documentclass[11pt, a4paper, twoside]{article}
\usepackage[a4paper]{anysize}
\usepackage{amsmath, amsfonts, amssymb, amsthm, graphicx}
\usepackage{amsmath,amssymb}
\usepackage[latin1]{inputenc}

\usepackage{amssymb}
\usepackage{amsmath}
\usepackage{graphics}
\usepackage{latexsym}
\usepackage{amsfonts}
\usepackage{color}
\usepackage{epstopdf}
\usepackage{ctable}
\usepackage{graphicx}
\usepackage{verbatim}
\newtheorem{theo}{Theorem}[section]

\newtheorem{cor}[theo]{Corollary}
\newtheorem{prop}[theo]{Proposition}
\newtheorem{uppg}[theo]{Exercise}
\newtheorem{remark}[theo]{Remark}
\newcommand{\be}{\begin{eqnarray*}}
\newcommand{\ee}{\end{eqnarray*}}
\newcommand{\ben}{\begin{eqnarray}}
\newcommand{\een}{\end{eqnarray}}

\def\subsecn (#1) {\medskip\ \ \ {\it #1}\medskip}

\newcommand{\lp}[1]{\left(\begin{array}{#1}}
\newcommand{\rp}{\end{array}\right)}

\newcommand{\leftd}[1]{\left\{\begin{array}{#1}}
\newcommand{\rightd}{\end{array}\right.}

\def\A {\mathbf{A}}

\def\C {\mathbf{C}}
\def\D {\mathbf{D}}

\def\F {\mathbf{F}}

\def\H {\mathbf{H}}
\def\I {\mathbf{I}}

\def\L {\mathbf{L}}
\def\M {\mathbf{M}}

\def\P {\mathbf{P}}
\def\Q {\mathbf{Q}}
\def\R {\mathbf{R}}
\def\S {\mathbf{S}}

\def\U {\mathbf{U}}
\def\V {\mathbf{V}}

\def\a {\boldsymbol{a}}

\def\c {\boldsymbol{c}}

\def\e {\boldsymbol{e}}

\def\l {\boldsymbol{l}}

\def\p {\boldsymbol{p}}
\def\q {\boldsymbol{q}}

\def\s {\boldsymbol{s}}

\def\u {\boldsymbol{u}}
\def\v {\boldsymbol{v}}
\def\w {\boldsymbol{w}}

\def\y {\boldsymbol{y}}

\def\sig {\mathbf{\Sigma}}

\def\Rb {\mathbb{R}}
\def\Pb {\mathbb{P}}
\def\Eb {\mathbb{E}}

\date{\today}
\begin{document}
\title{Estimation of noisy cubic spline using a natural basis\\
\small Azzouz Dermoune, Cristian Preda\\
Laboratoire Paul Painlev\'e, USTL-UMR-CNRS 8524.\\
UFR de Math\'ematiques, B\^at. M2, 59655 Villeneuve d'Ascq C\'edex, France.\\
azzouz.dermoune@univ-lille1.fr, cristian.preda@univ-lille1.fr}

\maketitle
{\bf Abstract.} 
We define a new basis of cubic splines such that the coordinates of a natural cubic spline 
are sparse. We use it to analyse and to extend the classical Schoenberg and Reinsch result
and to estimate a noisy cubic spline. We also discuss the choice of the smoothing parameter. All our results are illustrated graphically.   

{\bf keyword.} 
 Cubic spline. BLUP. Smoothing parameter. Natural cubic spline. Linear regression. 
 Bayesian model. Stein's Unbiased Risk Estimate. Prediction errors. 


\section{Introduction} 
We consider, for $n\geq 1$, the regression model 
\ben\label{infini} 
y_i=f(t_i)+w_i,\quad i=1, \ldots, n+1,
\een 
where $y_1, \ldots, y_{n+1}, t_1< \ldots < t_{n+1}$ are real-valued observations, $w_1, \ldots, w_{n+1}$ are measurement errors and 
$f:[t_1, t_{n+1}]\to \Rb$ is an unknown element of the infinite dimensional space $H^2$ of all functions with square integrable second derivative.
The approximation of $f$ by cubic splines considers the regression model 
\ben\label{fini} 
y_i=s(t_i)+w_i,\quad i=1, \ldots, n+1,
\een  
where $s$ is an unknown element of the finite dimensional space of cubic splines. 
Schoenberg \cite{schoenberg} introduced in 1946  the terminology spline for a certain type of piecewise polynomial interpolant. 
The ideas have their roots in the aircraft and shipbuilding industries. 
Since that time, splines have been shown to be applicable and effective for a large number of tasks in interpolation and approximation.
Various aspect of splines and their applications can be found in \cite{de Boor}, \cite{Claeskens}, \cite{Kramer}, \cite{Micula1999}, \cite{Paige} and \cite{reinsch}. 
See also the references therein.  

Let us first define properly the cubic splines approximation and introduce our notations.      
A map $s$ belongs to the set $S_3$ of cubic splines with the knots $t_1 <  . . . < t_{n+1}$ if there exist 
$(p_1, \ldots,p_{n+1})$ in $\Rb^{n+1}$, $(q_1, \ldots, q_n)$, $(u_1, \ldots, u_n)$, $(v_1, \ldots, v_n)$ in $\Rb^n$ such that, for $i=1, \ldots, n$ and $t\in [t_i,t_{i+1})$,  
\ben 
s(t)=p_i+q_i(t-t_i)+\frac{u_i}{2}(t-t_i)^2+\frac{v_i}{6}(t-t_i)^3. 
\label{cubicspline}
\een 
We are intereseted in the set $S_3\cap C^2$ of $C^2$-cubic splines. A cubic spline 
$s$, having its second derivatives $s''(t_1+)=s''(t_{n+1}-)=0$, is called natural.  
A well known result tells us that if $f\in H^2$ and $s\in S_3\cap C^2$ are such that $f(t_i)=s(t_i)$ for all 
$i=1, \ldots, n+1$, then $\int_{t_1}^{t_{n+1}}|s(t)-f(t)|^2dt=O(h^4)$ 
with $h=\max((t_{i+1}-t_i)^4:\, i=1, \ldots, n+1)$. See e.g. 
\cite{de Boor}, \cite{Sewell}. Hence, by paying the cost 
$O(h^4)$ we can replace the model (\ref{infini}) by (\ref{fini}).   

It is well known that any natural cubic spline of $S_3\cap C^2$ can be expressed using the 
{\it all} the  $n+3$ elements of the cubic B-spline basis, see e.g. \cite{Micula1999}. In Section 2 
we construct a new basis of $S_3\cap C^2$ in which 
any natural cubic spline needs only $n+1$ elements. In Sections 3-6 we treat the problem of 
estimation a noisy cubic spline.   
 
\section{The natural basis for $C^2$-Cubic splines}\label{sec2} 
Usually, the B-splines are used as a basis. The aim of this section is to construct a new basis 
which is more suitable for the natural cubic splines. Before going further, we need some notations. Let for $i=1, \ldots, n$, 
$h_i=t_{i+1}-t_i$.  
The spline $s$, defined in (\ref{cubicspline}), is of class $C^2$ if and only if   
\ben 
&&p_i+q_ih_i+\frac{u_i}{2}h_i^2+\frac{v_i}{6}h_i^3=p_{i+1},\quad i=1, \ldots, n,\label{czero}\\
&&q_i+u_ih_i+\frac{v_i}{2}h_i^2=q_{i+1}, \quad i=1, \ldots, n-1,\label{cone}\\
&&u_i+h_iv_i=u_{i+1}, \quad i=1, \ldots, n\label{ctwo}. 
\een 

We introduce the column vectors $\q=(q_1, \ldots,q_{n})^T$, 
$\p=(p_1, \ldots, p_{n+1})^T$, and $\u=(u_1, \ldots,u_{n+1})^T$, where 
$\M^T$ is the transpose of the matrix $\M$. Using (\ref{czero}), (\ref{cone}), (\ref{ctwo}), 
we can show that there exist 
three matrices $\Q, \U, \V$ such that 
\ben 
\q=\Q\lp{c} u_1\\ \p\\ u_{n+1}\rp,\label{q}\\
\u=\U\lp{c} u_1\\ \p\\ u_{n+1}\rp, \label{u}\\
\v=\V\lp{c} u_1\\ \p\\ u_{n+1}\rp.\label{v} 
\een 
See the appendices 1 and 2 for the details. 

Let us define, for each $i=1, \ldots, n$, the piecewise functions,  
\be 
&&\chi_i(t)=1_{[t_i,t_{i+1})}(t),\quad\chi_i^1(t)=(t-t_i)1_{[t_i,t_{i+1})}(t),\quad\chi_i^2(t)=(t-t_i)^21_{[t_i,t_{i+1})}(t),\\
&&\chi_i^3(t)=(t-t_i)^31_{[t_i,t_{i+1})}(t),\quad \chi_0=0,\quad\chi_{n+1}=1_{t_{n+1}},\quad\chi_{n+2}=0. 
\ee 
Here $1_A$ denotes the indicator function of the set $A$. 
Clearly, the set $[\chi_i, \chi_i^{k}: i=1, \ldots, n+1,\,k=1, 2, 3]$ forms a basis of the set of cubic splines  
$S_3$. 
The map $s$ has the coordinates $\p, \q, \u, \v$ in this basis, i.e.   
\be 
s=\lp{cccc}
[\chi_1 \ldots \chi_n\, \chi_{n+1}]&[\chi_1^1 \ldots \chi_n^1]&\frac{1}{2}[\chi_1^2 \ldots \chi_n^2\, 0]&\frac{1}{6}[\chi_1^3 \ldots \chi_n^3]\rp \lp{c}\p\\ \q\\ \u\\ \v\rp.
\ee
If $s$ is $C^2$, then from (\ref{q}), (\ref{u}), (\ref{v}) (see Appendix 1), we have 
\be 
s=([\chi_0\,\chi_1 \ldots \chi_n\, \chi_{n+1}\,\chi_{n+2}]+[\chi_1^1 \ldots \chi_n^1]\Q+\frac{1}{2}[\chi_1^2 \ldots \chi_n^2\, 0]\U+\frac{1}{6}[\chi_1^3 \ldots \chi_n^3]\V)\lp{c} u_1\\ \p\\ u_{n+1}\rp.
\ee 
The $C^2$ cubic spline $s$ can be rewritten in the following new basis: 
\be 
s=[\varphi_0\ldots \varphi_{n+2}]\lp{c} u_1\\ \p\\ u_{n+1}\rp,
\ee 
where, for $j=0, \ldots, n+2$, 
\ben 
\varphi_{j}=\chi_{j}+[\chi_1^1 \ldots \chi_n^1]\q_{\cdot j+1}+\frac{1}{2}[\chi_1^2 \ldots \chi_n^2\, 0]\u_{\cdot j+1}+\frac{1}{6}[\chi_1^3 \ldots \chi_n^3]\v_{\cdot j+1}. 
\label{newbasis}
\een 
Here $\a_{\cdot j}$ denotes the $j$th column of the matrix $\A$. Each element of the new basis is a $C^2$ cubic spline.

From (\ref{ui}), we derive that the set of natural cubic splines is spanned by the basis 
$(\varphi_{j}: j=1, \ldots, n+1)$.\\
$\bullet$ The spline $\varphi_0$ is the unique $C^2$ cubic spline interpolating the points  $(t_1,0)$, $\ldots$, $(t_{n+1},0)$ and 
such that $\varphi_0''(t_1+)=1$, $\varphi_0''(t_{n+1}-)=0$. Hence, $\varphi_0$ is not a natural cubic spline.\\
$\bullet$ The spline $\varphi_{j}$, for $j=1, \ldots, n+1$, is the unique natural cubic spline interpolating the points  
$(t_j,1)$, $((t_i,0), i\neq j)$.\\
$\bullet$ The spline $\varphi_{n+2}$ is the unique $C^2$ cubic spline interpolating the points  $(t_1,0)$, $\ldots$, $(t_{n+1},0)$ and 
such that $\varphi_{n+2}''(t_1+)=0$, $\varphi_{n+2}''(t_{n+1}-)=1$.  Hence, $\varphi_{n+2}$ is not a natural cubic spline.
 
\begin{figure}[!ht]
  \centering
     \includegraphics[width=8cm, height=12cm, angle=270]{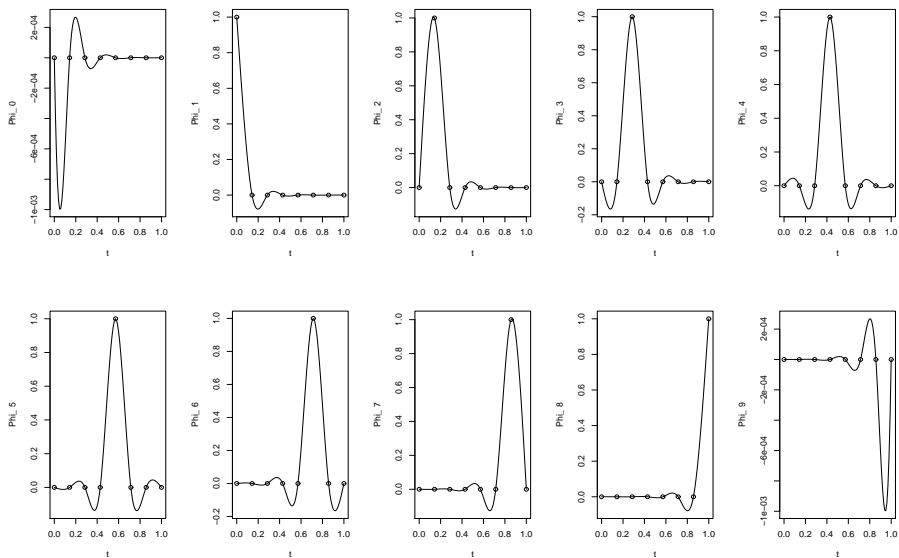}
 \caption{\small
The graph of the 10 elements of the natural basis. Here $n=7$ and $t_i=\frac{i-1}{n}$ for 
$i=1, \ldots,8$.}
\label{fig1}
\end{figure}

Observe that the natural cubic spline interpolating the points $(t_i,0),\,i=1, \ldots, n+1$ is the null map
\be 
s_0=[\varphi_0,\ldots, \varphi_{n+2}]\lp{c} 0\\ 0\\ 0\rp.
\ee 
As an illustration, in Figure \ref{fig1} we plot, for $n=7$, $t_i=\frac{i-1}{n}, i=1, \ldots, n+1$, the basis $\{\varphi_0, \ldots, \varphi_{n+2}\}$ and their derivatives in 
Figure \ref{fig2} and  Figure \ref{fig3}.
We can show that our basis has the reverse time property (see Figure \ref{fig1}), i.e. 
\be 
\varphi_j(t_{n+1}-t)=\varphi_{n+3-j}(t),\quad \forall\,j=0, \ldots, n+2,\quad t\in [t_1,t_{n+1}].
\ee 
Observe that our new basis is very different of the classical cubic B-spline basis. 
\begin{figure}[!ht]
\label{fig2}
  \centering
     \includegraphics[width=8cm, height=12cm, angle=270]{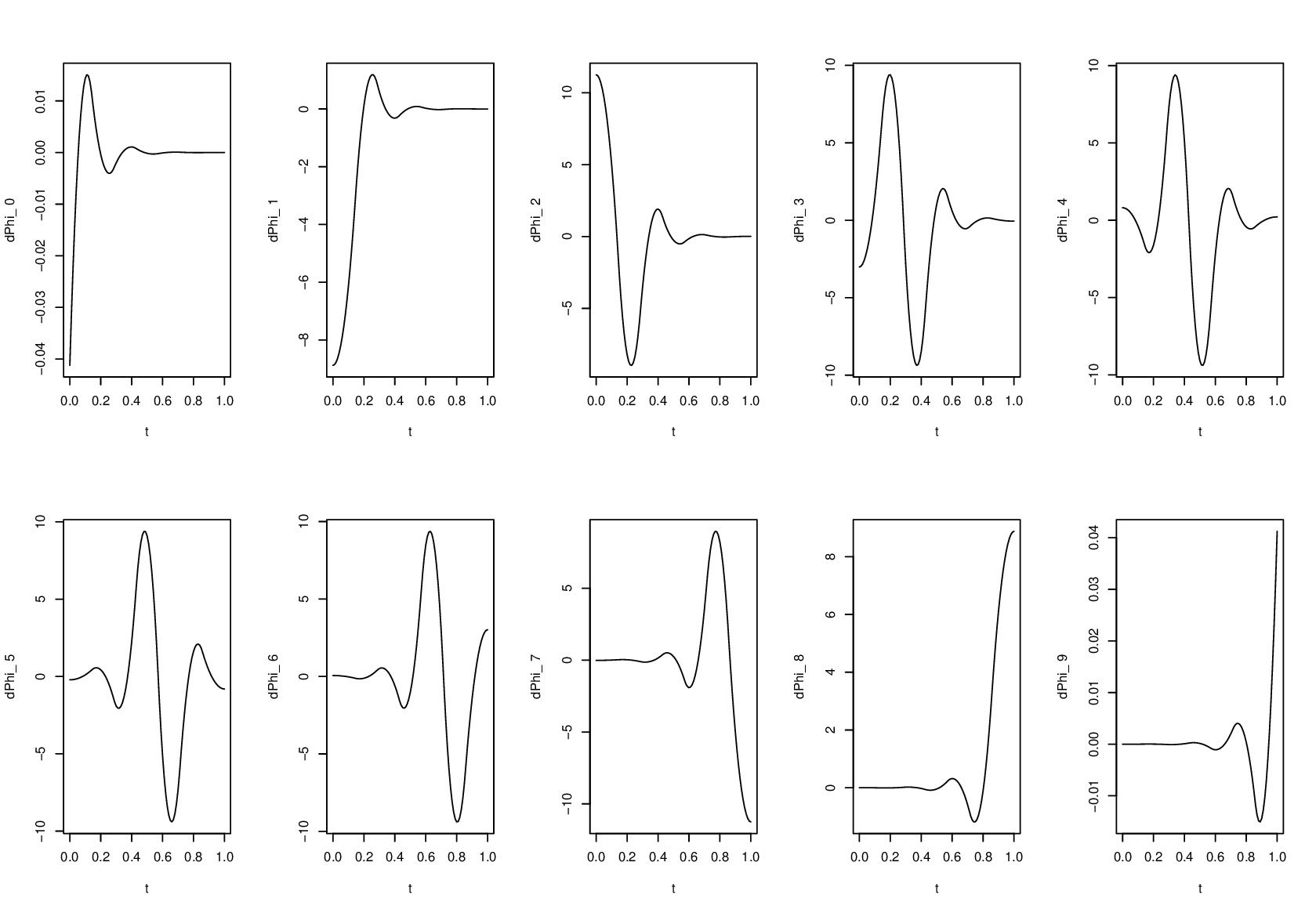}
 \caption{\small
The graph of the first derivative of the 10 elements of the natural basis. Here $n=7$ and $t_i=\frac{i-1}{n}$ for 
$i=1, \ldots,8$.}
\label{fig2}
\end{figure}
\begin{figure}[!ht]
  \centering
     \includegraphics[width=8cm, height=12cm, angle=270]{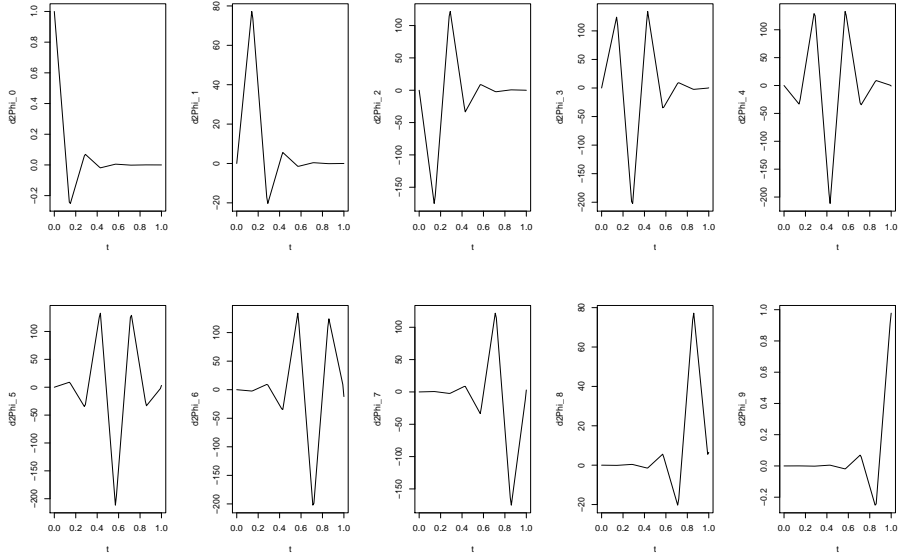}
 \caption{\small
The graph of the second derivative of the 10 elements of the natural basis. Here $n=7$ and $t_i=\frac{i-1}{n}$ for 
$i=1, \ldots,8$.}
\label{fig3}
\end{figure}

\section{Our basis and Schoenberg-Reinsch optimization}\label{interpolation}
In this section we use our new basis to review the well known results concerning the $L^2$ penalty and the optimal property of cubic splines. 
  
Let $p_i$, $i=1, \ldots, n+1,$ be a set of points in $\Rb$.   
The famous result of Schoenberg 1964 \cite{schoenberg} and Reinsch 1967 \cite{reinsch} tells us that the minimizer 
\be
I_2(\p)&\stackrel{def}{=}&\arg\min_{f\in H^2}\{\int_{t_1}^{t_{n+1}}|f''(t)|^2dt:\quad f(t_i)=p_i\quad, i=1, . . ., n+1\}\\
&=&\sum_{j=1}^{n+1}p_j\varphi_{j}
\ee 
is the natural $C^2$ cubic spline which interpolates the points $(t_i,p_i), i=1, \ldots, n+~1$. It follows, for $j=1, \ldots, n+1$, that  
\be 
I_2(\delta_j)&=&\arg\min_{f\in H^2}\{\int_{t_1}^{t_{n+1}}|f''(t)|^2dt:\quad f(t_i)=\delta_j^i\quad, i=1, . . ., n+1\}\\
&=&\varphi_j, 
\ee 
where $\delta_j\in\Rb^{n+1}$ and has the component $\delta_j^i=1$ if $i=j$ and $0$ otherwise. 

The aim of this section is to interpret Schoenberg and Reinsch result
using the natural basis. As a by-product, we will show that $\varphi_0$ and $\varphi_{n+2}$ are respectively solution of the following optimization problems:   
\ben 
\min_{f\in S_3\cap C^2}\{\int_{t_1}^{t_{n+1}}|f''(t)|^2dt:\quad f''(t_1)=1, \,\,f(t_i)=0\quad, i=1, . . ., n+1\},\\
\min_{f\in S_3\cap C^2}\{\int_{t_1}^{t_{n+1}}|f''(t)|^2dt:\quad f''(t_{n+1})=1, \,\,f(t_i)=0\quad, i=1, . . ., n+1\}.
\label{firstendsecondderivative} 
\een     

\subsection{Revisiting Schoenberg and Reinsch result} 

\begin{prop}\label{J2u}  
Let us introduce, for $\u\in\Rb^{n+1}$, the quadratic form  
\be 
J_2(\u)=\sum_{i=1}^{n}\frac{h_i}{3}(u_i^2+u_iu_{i+1}+u_{i+1}^2). 
\ee 
The minimization  
\be 
\min_{s\in C^2\cap S_3}\{\int_{t_1}^{t_{n+1}}|s''(t)|^2dt:\quad s(t_i)=p_i\quad, i=1, . . ., n+1\}
\ee
is equivalent to 
\be 
\min_{u_1, u_{n+1}}\{J_2(\U(u_1,\p,u_{n+1})^T)\}. 
\ee
\end{prop}
 
{\bf Proof.}   
Schoenberg and Reinsch result tells us that 
$I_2(\p)=\sum_{j=1}^{n+1}p_j\varphi_{j}$ is the minimizer of 
\ben 
\min_{s\in C^2\cap S_3}\{\int_{t_1}^{t_{n+1}}|s''(t)|^2dt:\quad s(t_i)=p_i\quad, i=1, . . ., n+1\}.
\label{two}
\een  
If $s\in C^2\cap S_3$, using (\ref{cubicspline}) and (\ref{vu}), then 
\be 
\int_{t_1}^{t_{n+1}}|s''(t)|^2 dt&=&\sum_{i=1}^{n}h_i\int_{0}^{1}|tu_i+(1-t)u_{i+1}|^2dt\\
&=&J_2(\u).
\ee
Now the equality (\ref{u}) achieves the proof. 

\subsection{Some consequences of Schoenberg-Reinsch result}
First, let us rewrite 
\ben\label{J2C} 
J_2(\u)&=&\sum_{i=1}^{n}\frac{h_i}{3}(u_i^2+u_iu_{i+1}+u_{i+1}^2)\\
&=&(u_1, \p^T, u_{n+1})\C\lp{c}u_1\\ \p \\u_{n+1}\rp, 
\een 
where 
\ben 
&&\C=\frac{h_1}{3}\U_1^T\U_1+\frac{2}{3}\sum_{i=2}^nh_i\U_i^T\U_i+\frac{h_{n}}{3}\U_{n+1}^T\U_{n+1}+\nonumber\\
&&\frac{1}{6}\sum_{i=1}^nh_i[\U_i^T\U_{i+1}+\U_{i+1}^T\U_i].
\label{Cmatrix}
\een 

Now, we summarize the properties of the matrix $\C$. 
\begin{prop}\label{matrixC} 
The matrix $\C$ is symmetric, and non-negative definite.  
The quadratic form $(u_1, \p^T, u_{n+1})\C\lp{c}u_1\\ \p \\u_{n+1}\rp=0$ 
if and only if $u_1=u_{n+1}=0$ and $\p$ belongs to the range $R(\L)$ of the matrix 
\ben\label{Lm} \L= 
\lp{cc} 1& t_1\\
\vdots&\vdots\\
1&t_{n+1}\rp.
\een  
It follows that, for all $j=1, \ldots, n+1$, that  $c_{1,j+1}=c_{n+3,j+1}=0$, i.e. the matrix  
\be 
\C=\lp{ccccc} c_{1,1}&0      &c_{1,n+3}\\ 
               0&\C(2,n+2)&0 \\
               c_{n+3,1}&0      &c_{n+3,n+3}
          \rp,
\ee  
where $\C(2,n+2)=[c_{ij}:i,j=2, \ldots n+2]$. The sub-matrix 
\ben\label{u1un+1} 
\lp{cc} c_{1,1}&c_{1,n+3}\\ 
      c_{n+3,1}& c_{n+3,n+3}\rp
\een 
is symmetric, positive definite. The null-space of the sub-matrix $\C(2,n+2)$ 
is equal to $R(\L)$. Moreover, from the decomposition 
$s=u_1\varphi_0+\sum_{j=1}^{n+1}p_j\varphi_{j}+u_{n+1}\varphi_{n+2}$, we have 
\ben 
\int_{t_1}^{t_{n+1}}|s''(t)|^2dt=u_1^2c_{1,1}+u_{n+1}^2c_{n+3,n+3}+2c_{1,n+3}u_1u_{n+1}+ \p^T\C(2,n)\p.
\label{statisticalinference}
\een 
\end{prop} 

From (\ref{statisticalinference}), we derive that the second derivatives $\{\varphi_j'': j=0, \ldots, n+2\}$ of the new basis 
satisfy 
\be 
\int_{t_1}^{t_{n+1}}\varphi_i''(t)\varphi_j''(t)dt=c_{i+1,j+1},\quad i,j=0,  \ldots, n+2. 
\ee 
We can show numerically that $\varphi_0''$ (respectively $\varphi_{n+2}''$) is orthogonal to $\varphi''_{j}$ for all $j=1, \ldots, n+2$ (respectively 
to $\varphi''_{j}$ for all $j=0, \ldots, n+1$). As an example    
the matrix $\C$ has the following form, for $n=7$, $t_i=\frac{i-1}{n}, i=1, \ldots, n+1$,
{\tiny
\ben\label{C}
\C = \left(\begin{array}{rrrrrrrrrr}
0.04 &    0.00&     0.00 &    0.00&     0.00 &    0.00 &    0.00 &    0.00 &    0.00 & 0.00 \\
\cline{2-9} 0.00  & 551.44& -1250.64 &  886.54 & -237.54&    63.63&   -16.97 &    4.24 &   -0.71 & 0.00\\
0.00& -1250.64 & 3387.82& -3261.27 & 1425.26 & -381.77&   101.80&   -25.45&     4.24&  0.00\\
 0.00 &  886.54& -3261.27&  4813.08& -3643.03&  1527.06&  -407.22 &  101.80 &  -16.97&  0.00\\
 0.00 & -237.54 & 1425.26& -3643.03&  4914.88& -3668.49&  1527.06&  -381.77&    63.63&  0.00\\
 0.00 &   63.63 & -381.77 & 1527.06& -3668.49&  4914.88& -3643.03&  1425.26&  -237.54&  0.00\\
 0.00 &  -16.97&   101.80&  -407.22&  1527.06& -3643.03&  4813.08 &-3261.27&   886.54&  0.00\\
 0.00&     4.24 &  -25.45&   101.80&  -381.77&  1425.26& -3261.27&  3387.82& -1250.64&  0.00\\
 0.00&    -0.71&     4.24&   -16.97&    63.63&  -237.54&   886.54& -1250.64&   551.44&  0.00\\  \cline{2-9}
0.00&     0.00&     0.00&     0.00&     0.00&     0.00&     0.00&     0.00&     0.00&  0.04\\
\end{array}\right).
\een 
}
Observe that the fact that sub-matrix (\ref{u1un+1}) is diagonal is not expected.  
 
Now we derive easily the following results. 
\begin{prop}\label{splinewithu1un}
 The splines $\varphi_0$ and $\varphi_{n+2}$ are respectively solution of the optimizations 
\ben
\min_{s\in S_3\cap C^2}\{\int_{t_1}^{t_{n+1}}|s''(t)|^2dt:\quad s''(t_1)=1, \,\,s(t_i)=0\quad, i=1, . . ., n+1\}\label{SC1}\\
\min_{s\in S_3\cap C^2}\{\int_{t_1}^{t_{n+1}}|s''(t)|^2dt:\quad s''(t_{n+1})=1, \,\,s(t_i)=0\quad, i=1, . . ., n+1\}\label{SCn+1} . 
\een
\end{prop} 

{\bf Proof.} The optimizations (\ref{SC1}),(\ref{SCn+1}) are equivalent to 
\be 
\min_{u_1,\p,u_{n+1}}\{(u_1, \p^T, u_{n+1})\C\lp{c}u_1\\ \p \\u_{n+1}\rp:\quad u_1=1, \quad \p=0\},\\
\min_{u_1,\p,u_{n+1}}\{(u_1, \p^T, u_{n+1})\C\lp{c}u_1\\ \p \\u_{n+1}\rp:\quad u_{n+1}=1, \quad \p=0\},
\ee 
and then have respectively the solutions, $(u_1=1, \p=0, u_{n+1}=0)$ and $(u_1=0, \p=0, u_{n+1}=1)$.  

More generaly we have the following result. 
\begin{prop}\label{splineparametrization}
 Let $k$ be a positive integer, $\M$ a $n+3$ by $n+3$ matrix, 
 $\A$ be a $k$ by $n+3$ matrix and $\lp{c}c_1\\ \vdots\\ c_{k}\rp:=\c\in \Rb^k$ all are given. 
Suppose that the null spaces $N(\A)$, $N(\M)$ do not overlap, i.e. 
$N(\A)\cap N(\M)=\{0\}$. 
Then the optimization 
\ben
\min_{u_1,\p,u_{n+1}}\{(u_1, \p^T, u_{n+1})\M\lp{c}u_1\\ \p \\u_{n+1}\rp:\quad\A\lp{c}u_1\\ \p \\u_{n+1}\rp=\c\}\label{GR}
\een
has a unique solution. More precisely, there exist a unique couple $(\l,\v)$
such that 
\be 
\M\v=\A^T\l,\quad \A\v=\c. 
\ee 
The vector $\v$ is the minimizer, and $\l$ is the Lagrange multiplier.  
\end{prop} 

\section{Natural cubic spline estimate} 
Let $s:\Rb\to \Rb$ be a natural cubic spline known with imprecision on the knots $t_1, \ldots, t_{n+1}$, i.e. 
\ben\label{noisynatural}
y_i=s(t_i)+w_i,\quad i=1, \ldots, n+1,
\een
where $w_i$ is the noise added to the true value $s(t_i)=p_i$. In the sequel $\y=(y_1, \ldots, y_{n+1})^T$.     
Schoenberg and Reinsch 
optimization, for each $\lambda > 0$,  
\be 
\arg\min_{f\in S_3\cap C^2}\{\lambda\int_{t_1}^{t_{n+1}}|f''(t)|^2dt+\sum_{i=1}^{n+1}|f(t_i)-y_i|^2\},
\ee 
provides an estimator of $s$. The parameter $\lambda> 0$ is called the smoothing parameter. 
Using the same arguments and notations as in Proposition \ref{J2u} and 
(\ref{J2C}) the latter optimization problem is equivalent to    
\ben
\min\{\lambda (u_1,\p^T,u_{n+1})\C\lp{c}u_1\\\p\\u_{n+1}\rp+\|\p-\y\|^2:\quad \p\in\Rb^{n+1}, u_1, u_{n+1}\in \Rb\}\label{naturalestimate}, 
\een
where $\|\cdot\|$ denotes the Euclidean norm. 

We have easily the following result.
\begin{prop}   
The equality  
\be 
(u_1,\p^T,u_{n+1})\C\lp{c}u_1\\\p\\u_{n+1}\rp=&&
u_1^2c_{11}+u_{n+1}^2c_{n+3,n+3}+2u_1u_{n+1}c_{1,n+3}+\\
&&\p^T\C(2,n+2)\p
\ee 
implies that the minimizer of (\ref{naturalestimate}) is $u_1=u_{n+1}=0$ and $\p$ is solution of the following system: 
\be 
(\I+\lambda\C(2,n+2))\p=\y.
\ee 
The solution $\p$ is given by   
\ben 
\hat{\p}=\H(\lambda)\y, 
\label{hatmatrix}
\een
where $\H(\lambda):=(\I+\lambda\C(2,n+2))^{-1}$ is called the hat matrix.
\end{prop} 

Now, we discuss the limits of (\ref{naturalestimate}) as 
$\lambda\to 0$ and $\lambda\to +\infty$. 
\begin{cor} 1) The problem (\ref{naturalestimate}), when $\lambda\to 0$, becomes 
\ben\label{interpolation} 
\min\{(u_1,\y^T,u_{n+1})\C\lp{c}u_1\\\y\\u_{n+1}\rp:\quad \p=\y, u_1, u_{n+1}\in \Rb\}. 
\een
Its minimizer is $u_1=u_{n+1}=0, \p=\y$ i.e.   
\be 
\lim_{\lambda\to 0+}\H(\lambda)=\I_{n+1}.
\ee 
where $\I_{n+1}$ is the $(n+1)\times (n+1)$ identity matrix. 

2) The problem (\ref{naturalestimate}), when $\lambda\to +\infty$, becomes 
\ben\label{Lregression} 
\min\{\|\y-\p\|^2: \quad 
\C\lp{c}u_1\\\p\\u_{n+1}\rp=0\}. 
\een
Its minimizer is $u_1=u_{n+1}=0$ and    
\be 
\p&=&\lim_{\lambda\to +\infty}\H(\lambda)\y\\
&=&\L(\L^T\L)^{-1}\L^T\y\\
&=&\pi_{R(\L)}\y, 
\ee 
where $\L$ is the linear model matrix (\ref{Lm}), and 
$\pi_{R(\L)}$ is the orthogonal projection on the range of $\L$. 
\end{cor} 

\begin{remark} 
We have easily $\C(2,n+2)\L=0$. From that we derive that $\H(\lambda)\L=\L$ and 
then $\L(\L^T\L)^{-1}\L^T[\H(\lambda)-\L(\L^T\L)^{-1}\L^T]=0$. It follows that 
$\L(\L^T\L)^{-1}\L^T\y$ and $[\H(\lambda)-\L(\L^T\L)^{-1}\L^T]\y$ are orthogonal. 
The component $Lreg\y:=\L(\L^T\L)^{-1}\L^T\y$ is the linear regression i.e. $Lreg\y$ is the orthogonal projection 
of the data $\y$ on the linear space $R(\L)$ (the range of $\L$). 
By introducing the orthogonal projections $\pi_{R(\L)}$ and 
$\pi_{R(\L)^{\perp}}$ respectively on $R(\L)$ and $R(\L)^{\perp}$, the minimizator  
\be
\H(\lambda)\y=\arg\min\{\lambda \p^T\C(2,n+2)\p+\|\p-\y\|_2^2:\quad \p\in\Rb^{n+1}\}
\ee 
is the sum of 
\be
\arg\min\{\lambda \p_1^T\C(2,n+2)\p_1+\|\p_1-\pi_{R(\L)^{\perp}}\y\|_2^2:\quad \p_1\in
R(\L)^{\perp}\},
\ee 
and 
\be
\arg\min\{\|\p_1-\pi_{R(\L)}\y\|^2:\quad \p_1\in R(\L)\}=\pi_{R(\L)}\y. 
\ee 
Hence, the component 
\be 
[\H(\lambda)-\L(\L^T\L)^{-1}\L^T]\y=\arg\min\{&&\lambda \p_1^T\C(2,n+2)\p_1+\|\p_1-\pi_{R(L)^{\perp}}\y\|_2^2:\\
&& \p_1\in R(L)^{\perp}\}
\ee 
is the penalized projection of $\y$ on $R(L)^{\perp}$, i.e.  
$[\H(\lambda)-\L(\L^T\L)^{-1}\L^T]\y$ is the nearest vector $\p\in R(L)^{\perp}$ 
to $\pi_{R(\L)^{\perp}}\y$ under the constraint $\p^T\C(2,n+2)\p\leq \delta$. Thanks 
to Lagrange multiplier, the positive 
constant $\delta$ and the smoothing parameter $\lambda$
are related by the equation $\y^T\H(\lambda)\C(2,n+2)\H(\lambda)\y=\delta$.                                         
The penalty $\p^T\C(2,n+2)\p=\|\C^{1/2}(2,n+2)\p\|^2$ measures the deviation 
of the vector $\p$ with respect to the linear space $R(\L)$.
The vector $\C^{1/2}(2,n+2)\p$ can be seen as an oblique projection 
on $R(\L)^{\perp}$. 
\end{remark} 
\section{Choice of the smoothing parameter $\lambda$} 
\subsection{Deterministic noise} 
We have, for any $\lambda >0$, the estimated model of (\ref{noisynatural}) 
\be 
\y=\H(\lambda)\y+[\I-\H(\lambda)]\y.   
\ee 
We proposed $\H(\lambda)\y$ as an estimator of $\p$ and therefore, 
$[\I-\H(\lambda)]\y$ is an estimator of 
the noise $\w=(w_1, \ldots,w_{n+1})^T$. 

The equality $[\I-\H(\lambda)]\y=\w$ holds only for $\lambda=0$, $\w=0$ and $\p$
is a straightline, i.e. $p_i=a+bt_i$ for all $i$. 
A natural way to link the smoothing parameter and the size of the noise is to solve the equation 
\ben  
\|\y-\H(\lambda)\y\|^2=\|\w\|^2.
\label{totalnoisesmoothing}
\een  
The following result shows that the equation (\ref{totalnoisesmoothing}) has a solution only for "small noise".     
\begin{prop} 
The map $\lambda\in (0, +\infty)\to \psi(\lambda)=\|\y-\H(\lambda)\y\|^2$ is concave, varies between $0$ and $\|\y-Lreg(\y)\|^2$.
The equation (\ref{totalnoisesmoothing}) has a solution if and only if 
\ben 
\|\y-Lreg(\y)\|^2> \|\w\|^2.
\label{total}
\een 
\end{prop} 

The proof is a consequence of the fact that $\H(\lambda)\to \I$ as $\lambda\to 0$ and  
$\H(\lambda)\to Lreg$ as $\lambda\to +\infty$.

\begin{figure}[!ht]
  \centering
     \includegraphics[width=8cm, height=12cm, angle=270]{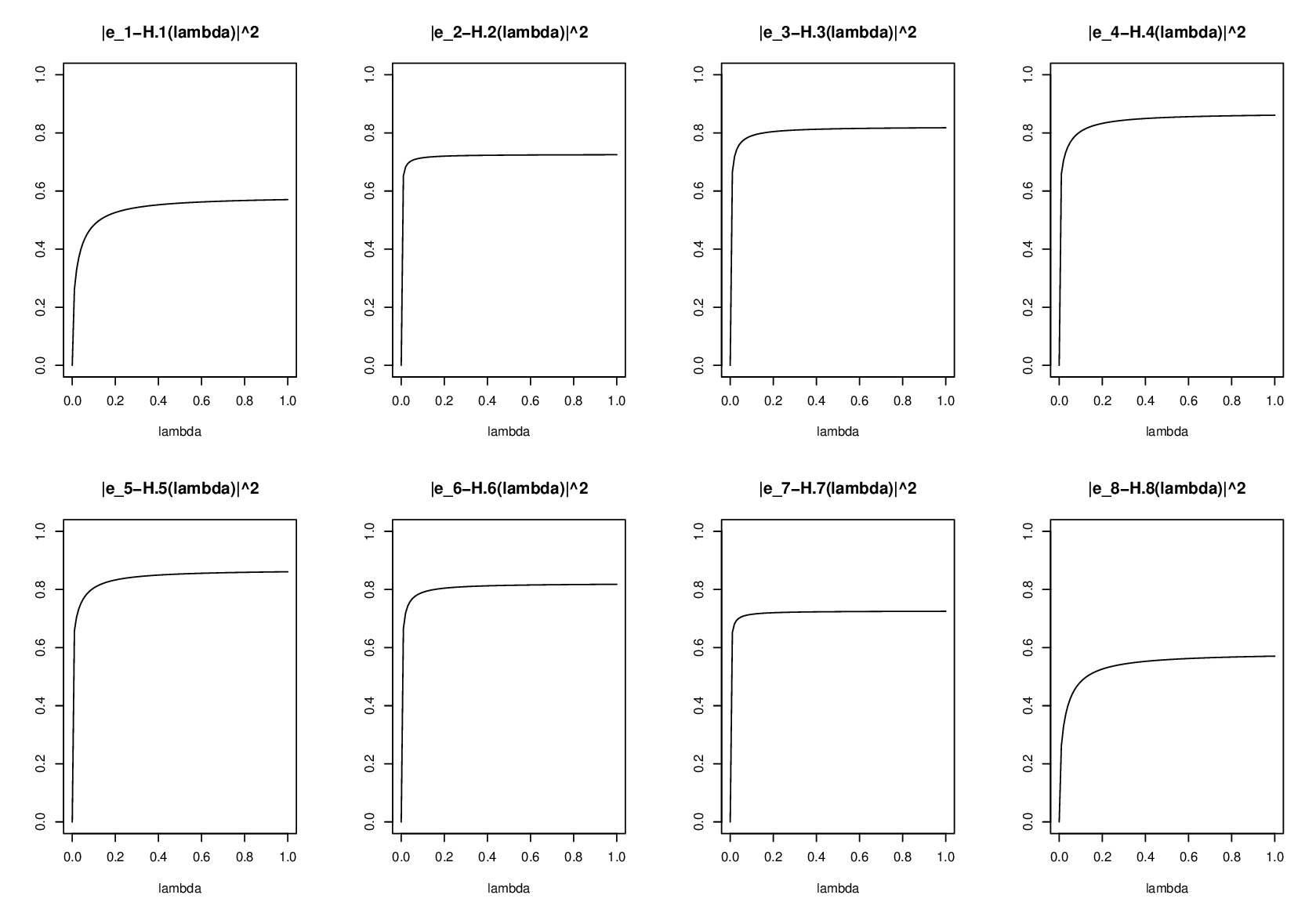}
 \caption{\small
 $\lambda\to \|\e_i-\H_{\cdot i}(\lambda)\|^2$.}
\label{fig4}
\end{figure}

Observe that (\ref{total}) is equivalent to 
\ben 
\sum_{i=1}^{n+1}y_i^2\|\e_i-Lreg(\e_i)\|^2+2\sum_{i< j}y_iy_j\langle\e_i-Lreg(\e_i),\e_j-Lreg(\e_j)\rangle > \|\w\|^2,
\label{ijSNR}
\een 
and the smoothing parameter is solution of the equation 
\ben 
\sum_{i=1}^{n+1}y_i^2\|\e_i-\H_{\cdot i}(\lambda)\|^2+2\sum_{i< j}y_iy_j\langle\e_i-\H_{\cdot i}(\lambda),\e_j-\H_{\cdot j}(\lambda)\rangle = \|\w\|^2.
\label{ijSNRlambda}
\een 
It follows that the size of the weights $\|\e_i-Lreg(\e_i)\|^2$, $\langle\e_i-Lreg(\e_i),\e_j-Lreg(\e_j)\rangle$, 
$\|\e_i-\H_{\cdot i}(\lambda)\|^2$, $\langle\e_i-\H_{\cdot i}(\lambda),\e_j-\H_{\cdot j}(\lambda)\rangle$, 
are crucial in the existence of the smoothing parameter (\ref{totalnoisesmoothing}). 
In Figure \ref{fig4} we plot for $i=1, \ldots,8$ the graph of $\lambda\to \|\e_i-\H_{\cdot i}(\lambda)\|^2$ for $n=7$ and $t_i=\frac{i-1}{n}$. 

Remark that for the model $\y=y_i\e_i$ the $i$th column $\H_{\cdot i}(\lambda)$ is an estimator of the signal. Thanks to the equation (\ref{ijSNRlambda}), the quantity $\|\e_i-\H_{\cdot i}(\lambda)\|^2$ represents the noise-to-signal ratio (NSR), i.e.  the smoothing parameter is solution of
\ben 
\|\e_i-\H_{\cdot i}(\lambda)\|^2=\frac{\|\w\|^2}{y_i^2}. 
\label{NSRi}
\een  
For large noise, there is no smoothing parameter solution of (\ref{NSRi}). 
 
\begin{figure}[!ht]
  \centering
     \includegraphics[width=8cm, height=12cm, angle=270]{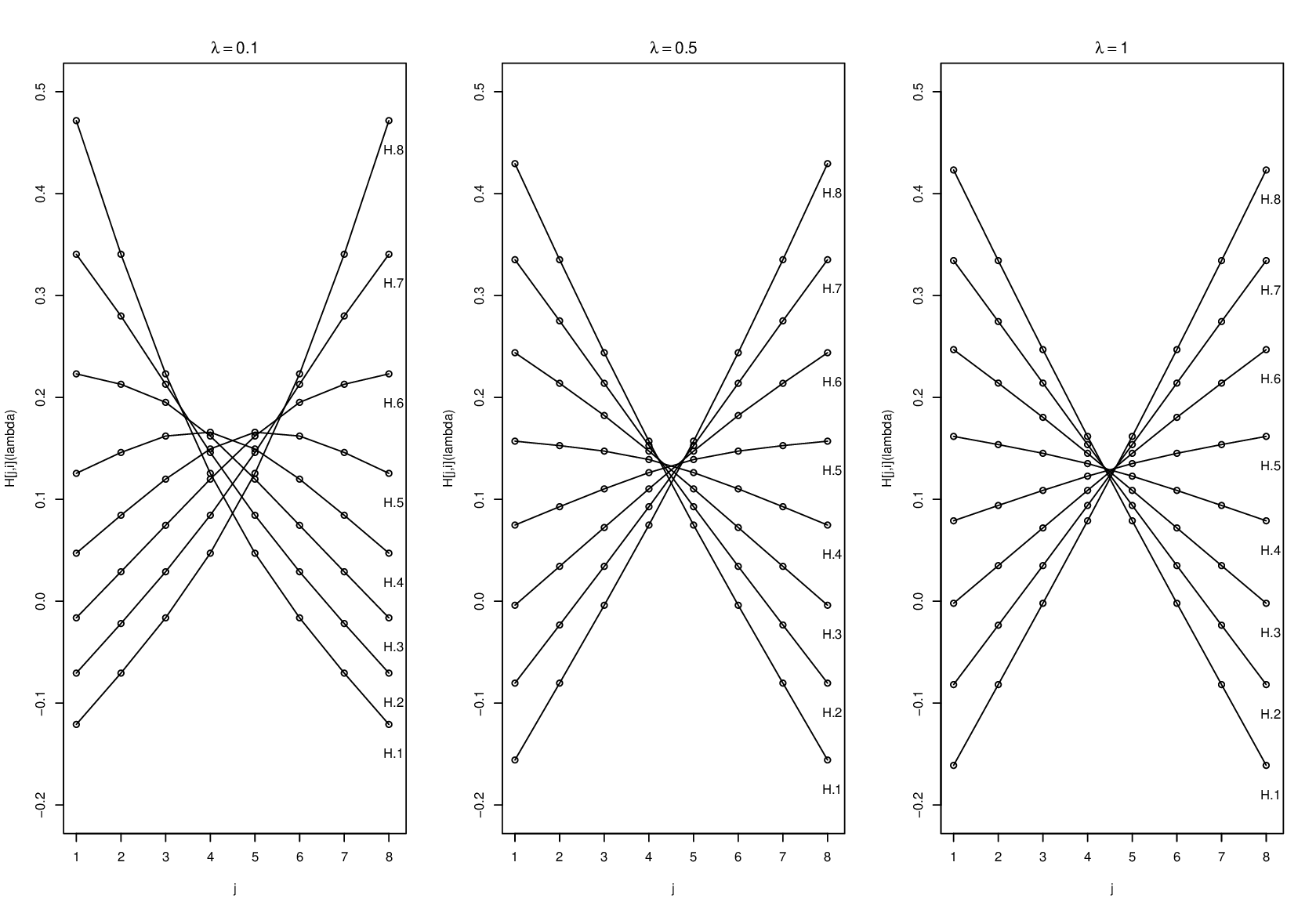}
 \caption{\small
Plot of $j\to \H_{ji}(\lambda)$. Here $n=7$ and $t_i=\frac{i-1}{n}$ for 
$i, j=1, \ldots,8$.}
\label{fig10}
\end{figure}

In Figure \ref{fig10} we plot the ``rainbow'' $\H_{\cdot i}(\lambda)$ for $i=1, \ldots, n+1$ and $\lambda\in \{0.1, 0.5, 1\}$.   

Concerning the weights $\|\e_i-Lreg(\e_i)\|^2$, $\langle\e_i-Lreg(\e_i),\e_j-Lreg(\e_j)\rangle$, we can calculate them explicitly as following. 
We recall that 
the linear regression of the data $\e_i$ is given by 
\be 
Lreg(\e_i)(t)= \beta_1^i+t\beta_2^i, 
\ee 
with $\beta_2^i=\frac{t_i-\bar{t}}{\sum_{j=1}^{n+1}(t_j-\bar{t})^2}$, 
$\beta_1^i=\frac{1}{n+1}-\beta_2^i\bar{t}$, and 
$\bar{t}$ denotes the empirical mean of the knots. 
The straightlines $(t\to Lreg(\e_i)(t): i=1, \ldots, n+1)$ have the common point $(\bar{t}, \frac{1}{n+1})$. 
Moreover we have, for $i\neq j$, that   
\be 
\|\e_i-Lreg(\e_i)\|^2&=&1-Lreg(\e_i)(t_i),\\
\langle\e_i-Lreg(\e_i),\e_j-Lreg(\e_j)\rangle&=&-Lreg(\e_i)(t_j),\\
&=&-Lreg(\e_j)(t_i).
\ee   
From all that we get the following result. 
\begin{prop}\label{RSSei} We have, for $i\neq j$,  
\be 
\|\e_i-Lreg(\e_i)\|^2 &=&\frac{n}{n+1}-\frac{(t_i-\bar{t})^2}{\sum_{j=1}^{n+1}(t_j-\bar{t})^2},\\
\langle\e_i-Lreg(\e_i),\e_j-Lreg(\e_j)\rangle &=&-\frac{1}{n+1}-\frac{(t_i-\bar{t})(t_j-\bar{t})}{\sum_{j=1}^{n+1}(t_j-\bar{t})^2}.
\ee 
\end{prop} 

It follows that the most important weights $\|\e_i-Lreg(\e_i)\|^2$ are when the $t_i$'s are close to $\bar{t}$. 
The most important negative correlation  
$\langle\e_i-Lreg(\e_i),\e_j-Lreg(\e_j)\rangle$ is given by the couple of end-points $(t_1,t_{n+1})$. 
The most important positive correlations  
$\langle\e_i-Lreg(\e_i),\e_j-Lreg(\e_j)\rangle$ are given by the begining $(t_1,t_2)$ and the ending $(t_{n},t_{n+1})$ of the knots.
The message of these remarks is that the allowed size of the noise depends on the values of data at the end-points $(t_1,t_{n+1})$ and at center i.e. near $\bar{t}$.

\begin{figure}[!ht]
  \centering
     \includegraphics[width=8cm, height=12cm, angle=270]{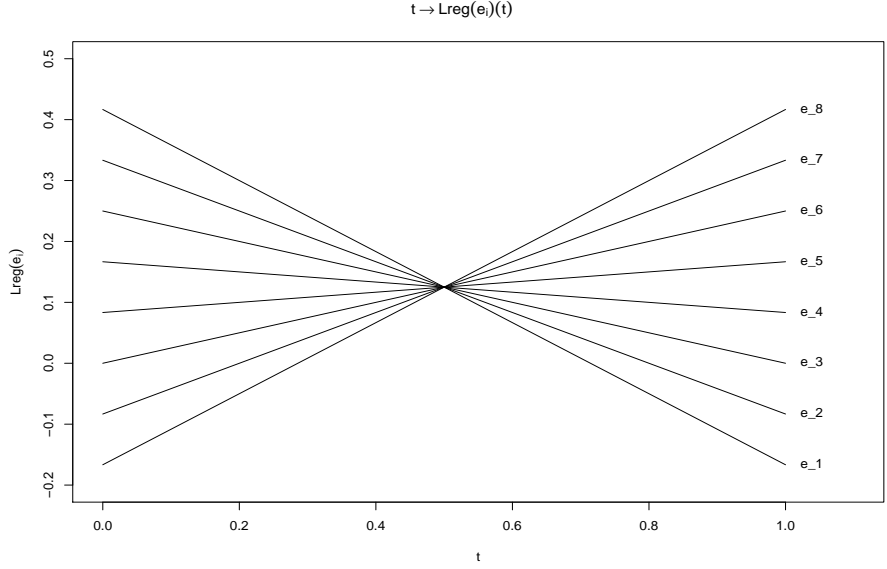}
 \caption{\small
Plot of $t\to Lreg(\e_i)(t)$. Here $n=7$ and $t_j=\frac{j-1}{n}$ for 
$i, j=1, \ldots,8$.}
\label{fig9}
\end{figure}

Figure \ref{fig9} shows that the straightlines $(t\to Lreg(\e_j)(t): j=1, \ldots, n+1)$ turn in the trigonometric sense around their common point. Remark that  the Figure \ref{fig10} 
illustrates also the convergence of $\H_{\cdot i}(\lambda)\to [Lreg(\e_i)(t_j): j=1, \ldots, n+1]^T$ as $\lambda\to +\infty$. 

What can we do if the condition (\ref{total}) does not hold ? In this case 
for all $\lambda\geq 0$, $\|\y-\H(\lambda)\y\|^2$ represents only a part of the noise i.e.  
\be 
\frac{\|\y-\H(\lambda)\y\|^2}{\|\w\|^2} < \frac{\|\y-Lreg\y\|^2}{\|\w\|^2}\in [0, 1).
\ee 
\subsection{Gaussian white noise} 
We suppose that the noise $\w$ is Gaussian and   
white with the variance $\sigma_w^2$. 
In this case $\sum_{i=1}^{n+1}\frac{w_i^2}{\sigma_w^2}$ has the $\chi_{n+1}^2$ probability distribution 
($(\frac{w_i}{\sigma}:\quad i=1, \ldots, n+1)$ are i.i.d. with the common distribution $\mathcal{N}(0,1)$,
the standard Gaussian distribution). 
For all $\varepsilon >0$, the event  
\be 
(n+1)(1-\varepsilon)\leq \frac{\|\w\|^2}{\sigma^2}\leq (n+1)(1+\varepsilon)
\ee 
holds with the probability 
\be 
\Pb((1-\varepsilon)\leq \frac{\chi_{n+1}^2}{n+1}\leq (1+\varepsilon)).
\ee 
The latter probability is close to 1 as $n$ becomes large. 

A first way to link the smoothing parameter to the noise is to choose $\lambda$ solution of the following constraint 
\ben
(n+1)(1-\varepsilon)\sigma^2\leq \|\y-\H(\lambda)\y\|^2\leq (n+1)(1+\varepsilon)\sigma^2.
\label{lambdavarepsilon} 
\een 
We denote respectively $\lambda^-(\sigma^2,\varepsilon,n+1)$, $\lambda^+(\sigma^2,\varepsilon,n+1)$
the solution of the equations 
\ben 
\|\y-\H(\lambda)\y\|^2=(n+1)(1-\varepsilon)\sigma^2,
\label{lambdaminus} 
\een 
and 
\ben 
\|\y-\H(\lambda)\y\|^2=(n+1)(1+\varepsilon)\sigma^2.
\label{lambdaplus} 
\een 
The solution of (\ref{lambdaminus}) exists under the hypothesis 
\ben 
(1-\varepsilon)\sigma^2 < \frac{\|\y-Lreg\y\|^2}{n+1}.
\label{minusvarianceconstraint}
\een 
The solution of (\ref{lambdaplus}) exists under the hypothesis 
\ben 
(1+\varepsilon)\sigma^2 < \frac{\|\y-Lreg\y\|^2}{n+1}.
\label{plusvarianceconstraint}
\een 
Remark that if $\lambda^+(\sigma^2,\varepsilon,n+1)$ exists then $\lambda^-(\sigma^2,\varepsilon,n+1)$ also exists. 
But in general the opposite is false. To understand the constraints (\ref{minusvarianceconstraint}) and (\ref{plusvarianceconstraint}),
we are going to study the quantity $\|\y-Lreg\y\|^2$ as a function of the signal $\p$ and the noise $\w$.

If the model is $y_i=a+bt_i+w_i$ with $i=1,\ldots, n+1$, then $Lreg\y=\H(+\infty)\y$ is 
the maximum likelihood estimator of the vector $\L(a,b)^T$.  
Moreover, $\frac{\|\y-Lreg\y\|^2}{n-1}$ is an unbiased consistent estimator of the variance $\sigma^2$. More precisely, $\frac{\|\y-Lreg\y\|^2}{\sigma^2}$ has  
the $\chi_{n-1}^2$-distribution. Hence,      
the constraint (\ref{plusvarianceconstraint}) holds with the probability 
$\Pb(\chi_{n-1}^2>(n+1)(1+\varepsilon))\to 0$ as $n\to +\infty$. 
But for $\varepsilon > \frac{2}{n+1}$, the constraint (\ref{minusvarianceconstraint}) holds with  
the probability 
$\Pb(\chi_{n-1}^2>(n-1)(1-\varepsilon))\to 1$ as $n\to +\infty$. 

In the general case we have the following result.  
\begin{prop}\label{noiseandsignal} 
Let $\y=\p+\w$, where $\w$ is the Gaussian white noise with the variance $\sigma_w^2$. We have  
\be 
\Eb[\|\y-Lreg\y\|^2]=(n-1)\sigma_w^2+\|\p-Lreg\p\|^2. 
\ee  
If the noise is fixed, then $\p\to \Eb[\|\y-Lreg\y\|^2]$ is minimal at the straightlines, i.e. $p_i=a+bt_i$
for all $i=1, \ldots, n+1$. The minimal value is equal to 
\ben 
\Eb[\|\y-Lreg\y\|^2]=(n-1)\sigma^2.
\label{unbiasednoiseestmator}
\een 
\end{prop} 

{\bf Proof.} From the equality $\y=\p+\w$, we have 
\be 
\|\y-Lreg\y\|^2=\|\p-Lreg\p\|^2+\|\w-Lreg \w\|^2+2\langle \p-Lreg\p,\w-Lreg\w\rangle.
\ee 
The rest of the proof is consequence of $\Eb(\w)=0$ and $\Eb[\|\w-Lreg\w\|^2]=(n-1)\sigma^2$. 

Roughly speaking,  Proposition \ref{noiseandsignal} combined with (\ref{minusvarianceconstraint}) and (\ref{plusvarianceconstraint})
tell us that the smoothing parameter 
\be
\lambda\in [\lambda^-(\sigma^2,\varepsilon,n+1),\lambda^+(\sigma^2,\varepsilon,n+1)]
\ee 
exists under the constraint 
\be 
\|\p-Lreg\p\|^2\approx 2 \sigma^2.
\ee 
\subsection{Smoothing parameter, SURE and PE} 
A second way to choose the smoothing parameter is to consider Stein's unbiased risk estimate (SURE) and the predictive risk error (PE).\\ 
a) Stein's Unbiased Risk estimate (SURE) \cite{Efron}, \cite{Stein}: The quadratic loss of the estimation of the vector $\p$ by $\H(\lambda)\y$ is equal to 
\be 
\|\H(\lambda)\y-\s\|^2=\sum_{i=1}^{n+1}|\H(\lambda)\y(i)-s(t_i)|^2,
\ee 
and the residual sum of squares is defined by 
\be 
RSS(\lambda):=\|\y-\H(\lambda)\y\|^2.
\ee 
The mean square risk is equal to 
\be 
&&R(\H(\lambda)\y,\p)=\Eb[\|\H(\lambda)\y-\p\|^2]\\
&&=\Eb[\|\y-\p\|^2]+\Eb[RSS(\lambda)]-2cov(\w-\H(\lambda)\w,\w)\\
&&=\Eb[\|\y-\p\|^2]+\Eb[RSS(\lambda)+2\sigma^2 (Trace(\H(\lambda))-(n+1))]\\
&&=\Eb[RSS(\lambda)+2\sigma^2Trace(\H(\lambda))-(n+1)\sigma^2]. 
\ee 
The quantity 
\be 
RSS(\lambda)+2\sigma^2Trace(\H(\lambda))-(n+1)\sigma^2
\ee 
is an unbiased risk estimate (called Stein's Unbiased Risk estimate, SURE for short). 
By minimizing SURE with respect to $\lambda\in (0,+\infty)$ we provide a criterion for choosing the  smoothing parameter
$\lambda_{SURE}$. 

b) Prediction and Training errors (PE). The prediction error is our error on a new observations 
$y_i^*=s(t_i)+w^*(t_i), i=1, \ldots, n+1$ independent of $\y$. If we predict the vector $\p$ by $\M(\lambda)\y$, then  the predictive     
risk PE is equal to 
\be 
&&\Eb[\|\y^*-\H(\lambda)\y\|^2]=\Eb[\|\p-\H(\lambda)\y\|^2]+(n+1)\sigma^2\\
&&=R(\H(\lambda)\y,\p)+(n+1)\sigma^2\\
&&=\Eb[RSS(\lambda)+2\sigma^2Trace(\H(\lambda))]. 
\ee 
Hence, $RSS(\lambda)+2\sigma^2Trace(\H(\lambda))$ is an unbiased estimate of the prediction error.
It follows that minimizing SURE is equivalent to minimize PE and then $\lambda_{SURE}=\lambda_{PE}$.  

\begin{figure}[!ht]
  \centering
     \includegraphics[width=8cm, height=12cm, angle=360]{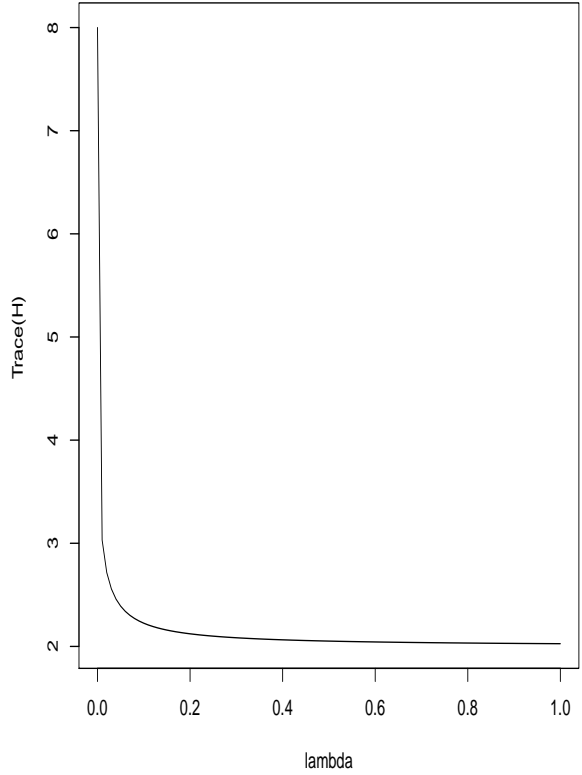}
 \caption{\small
Plot of $\lambda\in (0, +\infty)\to Trace(\H(\lambda))$.}
\label{fig11}
\end{figure}

In Figure\ref{fig11} we plot, for $n=7$, $i=1, \ldots, n+1$, $t_i=\frac{i-1}{n}$, the map 
$\lambda\in (0, +\infty)\to Trace(\H(\lambda))$. 

\section{Cubic spline estimate: General case} 
In this section we propose to find suitable symmetric and non-negative definite matrices 
$\P_{pen}=[p_{ij}: i,j=1, \ldots, n+3]$ such that the minimizer
\ben
(\hat{u}_1, \hat{p}_1, \ldots, \hat{p}_{n+1}, \hat{u}_{n+1})=\arg\min\{\lambda (u_1,\p^T,u_{n+1})\P_{pen}\lp{c}u_1\\ \p \\ u_{n+1}\rp+\nonumber\\
\|\y-\p\|^2:\quad u_1, p_1, \ldots, p_{n+1}, u_{n+1}\}\label{Plambda}
\een 
is a non natural cubic spline, i.e. $(\hat{u}_1,\hat{u}_{n+1})\neq (0,0)$. 
The following proposition addresses the uniqueness and the capacity  
of the estimator (\ref{Plambda}) to rediscover a non natural spline.   
  
\begin{prop} 
1) The minimizer of (\ref{Plambda}) is unique if and only 
if the sub-matrix $\lp{cc}p_{1,1}&p_{1,n+3}\\
p_{n+3,1}&p_{n+3,n+3}\rp$ is invertible. In this case 
the minimizer is given by 
\ben\label{Plambdaformula}
\lp{c}\hat{u}_1\\ \hat{\p} \\ \hat{u}_{n+1}\rp&=&
(\lambda \P_{pen}+ \Pi^T\Pi)^{-1}\lp{c}0\\ \y\\ 0\rp\\
:&=&\H_{\P_{pen}}(\lambda)\lp{c}0\\ \y\\ 0\rp 
\een 
where $\Pi\lp{c}u_1\\ \p\\ u_{n+1}\rp=\p$ for all $u_1, \p, u_{n+1}$.  

2) The condition 
\ben\label{u1}
(p_{1,j}:\quad j=2, \ldots, n+2)\neq 0,
\een 
respectively
\ben \label{un+1}   
(p_{n+3,j}:\quad j=2, \ldots, n+2)\neq 0, 
\een 
is the necessary condition 
which guaranties that $\hat{u}_1\neq 0$ respectively $\hat{u}_{n+1}\neq 0$.
\end{prop}  

Now, we discuss the limits of (\ref{Plambdaformula}) as 
$\lambda\to 0$ and $\lambda\to +\infty$. 
\begin{cor}\label{Pintreg} 
1) The limit 
\ben\label{Pinterpolation} 
\H_{\P_{pen}}(\lambda)\lp{c}0\\ \y\\ 0\rp\to 
\lp{c}u_1^0\\ \y\\ u_{n+1}^0\rp
\een 
as $\lambda \to 0$. Here $u_1^0, u_{n+1}^0$ is a minimizer 
of the objective function  
\be
(u_1,u_{n+1})\to (u_1,\y^T,u_{n+1})\P_{pen}\lp{c}u_1\\ \y \\ u_{n+1}\rp,
\ee
i.e. $u_1, u_{n+1}$ is solution of the linear system 
\be 
u_1p_{1,1}+u_{n+3}p_{1,n+3}=-\sum_{i=1}^{n+1}y_ip_{1,i+1}\\
u_1p_{n+3,1}+u_{n+3}p_{n+3,n+3}=-\sum_{i=1}^{n+1}y_ip_{n+3,i+1}. 
\ee 
In particular, if the data $\y$ is not orthogonal to the space spanned by 
the vectors $(p_{1,j}:\quad j=2, \ldots, n+2)^T$,  
$(p_{n+3,j}:\quad j=2, \ldots, n+2)^T$, then 
$\hat{u}_1$, $\hat{u}_{n+1}$ can't be both equal to zero. Namely, the estimator (\ref{Pinterpolation}) is not a natural spline. 

2) The limit 
\ben 
\H_{\P_{pen}}(\lambda)\lp{c}0\\ \y\\ 0\rp\to \arg\min\{\|\y-\p\|^2:\quad \P_{pen}\lp{c}u_1\\ \p \\ u_{n+1}\rp=0\}\label{Pregression}, 
\een 
as $\lambda \to +\infty$. 
\end{cor} 

{\bf Examples.} 
The matrix $\P_{pen}=\C$ given in (\ref{Cmatrix}) corresponds to the penalty 
\be 
(u_1,\p^T,u_{n+1})^T\C(u_1,\p^T,u_{n+1})^T=\int_{t_1}^{t_{n+1}}|s''(t)|^2dt
\ee 
where the cubic spline $s=u_1\varphi_0+\sum_{i=1}^{n+1}p_j\varphi_i+u_{n+1}\varphi_{n+2}$. 
A natural way to construct new matrices $\P_{pen}$ is to consider the more general penalization     
\be 
\int_{t_1}^{t_{n+1}}|a_2s''(t)+a_1s'(t)+a_0s(t)|^2dt=
(u_1,\p^T,u_{n+1})^T\P_{pen}(u_1,\p^T,u_{n+1})^T, 
\ee
where $a_0, a_1, a_2$ are given real numbers. 
Hence 
\ben\label{a}
\P_{pen}=a_0^2\C_{00}+a_1^2\C_{11}+a_2^2\C_{22}+a_0a_1[\C_{01}+\C_{01}^T]+a_0a_2[\C_{02}+\C_{02}^T]+a_1a_2[\C_{12}+\C_{12}^T], 
\een 
where $\C_{00}, \C_{01}, \C_{02}, \C_{11}, \C_{12}$ are respectively defined by 
\be 
\C_{00}&=&[\int_{t_1}^{t_{n+1}}\varphi_{i-1}(t)\varphi_{j-1}(t) dt:\quad i,j=1, \ldots, n+3],\\
\C_{01}&=&[\int_{t_1}^{t_{n+1}}\varphi_{i-1}(t)\varphi_{j-1}'(t) dt:\quad i,j=1, \ldots, n+3],\\
\C_{02}&=&[\int_{t_1}^{t_{n+1}}\varphi_{i-1}(t)\varphi_{j-1}''(t) dt:\quad i,j=1, \ldots, n+3],\\
\C_{11}&=&[\int_{t_1}^{t_{n+1}}\varphi_{i-1}'(t)\varphi_{j-1}'(t) dt:\quad i,j=1, \ldots, n+3],\\
\C_{12}&=&[\int_{t_1}^{t_{n+1}}\varphi_{i-1}'(t)\varphi_{j-1}''(t) dt:\quad i,j=1, \ldots, n+3],
\ee 
and the matrix $\C_{22}=\C$ defined in (\ref{Cmatrix}). 
 
The matrix $\P_{pen}$, for $a_0=a_1=a_2=1$, $n=7$, $t_i=\frac{i-1}{n}, i=1, \ldots, n+1$,
has the following form.

{\tiny
$$\mathbf{P}_{\mbox{pen}}=\left(\!\!\!\!\begin{array}{rrrrrrrrrr}
2.6&-239.22&614.68&-559.93&255.91&-96.85&33.253&-9.951&1.841&-0.01\\
-239.22&30883.64&-89524.08&98923.05&-57634.4&23760.7&-8518.95&2611.28&-488.75&1.84\\
614.68&-89524.08&278142.5&-345243.9&238016.7&-113691.9&43414.47&-13742.88&2611.71&-9.94\\
-559.93&98923.05&-345243.9&516426.7&-459007.5&281450.3&-127439.9&43417.03&-8521.09&33.24\\
255.91&-57634.4&238016.7&-459007.5&559860.4&-472755.5&281452.9&-113702.2&23768.8&-96.81\\
-96.85&23760.7&-113691.9&281450.3&-472755.5&559863&-459017.7&238055.1&-57664.7&255.76\\
33.253&-8518.95&43414.47&-127439.9&281452.9&-459017.7&516465.1&-345387.2&99036.12&-559.40\\
-9.951&2611.28&-13742.88&43417.03&-113702.2&238055.1&-345387.2&278677.6&-89946.04&612.72\\
1.841&-488.75&2611.71&-8521.09&23768.8&-57664.7&99036.12&-89946.04&31219.29&-236.99\\
-0.01&1.84&-9.95&33.24&-96.819&255.76&-559.40&612.722&-236.99&3.85
\end{array}\!\!\!\!\right)$$
}
 Observe that the conditions (\ref{u1}) and (\ref{un+1}) 
are satisfied. 
\section{Bayesian Model and statistical analysis}
The aim of this section is to give the Bayesian interpretation of the matrix 
penalization $\P_{pen}$. Let us first set the noisy cubic spline 
estimate in the context of the general linear model:   
\be 
\y=\F\beta+\R\eta+\w,
\ee   
where $\beta$ is an unknown parameters, $\F$ and $\R$ are known matrices. The random effects $\eta$ and the noise $\w$ are unknown, centred and independent random vectors. Their covariance matrices
 $cov(\eta):=\sig_{\eta}$, $cov(\w):=\sig_w$ are known. The term $\F\beta$ is called the fixed effects and $\R\eta$ is the random effects.

Let us revisit the 
best linear unbiased predictors (BLUP) and the best linear unbiased estimators
(BLUE).  
There is a long history and huge literature on this subject, see for instance
\cite{DRW}, \cite{DDR1}, \cite{DDR2}, \cite{Henderson1}, \cite{Henderson2}, \cite{Henderson3}, \cite{Henderson4}, \cite{Rao}, \cite{Robinson}, \cite{3chinois}
and references herein.

 {\bf BLUE of $\beta$.} The BLUE of $\beta$ is the estimator $\hat{\beta}=\hat{\M}_{\beta}\y$, 
with $\hat{\M}_{\beta}$ (called the hat matrix of $\hat{\beta}$) being the matrix such that $\hat{\M}_{\beta}\F=\I$ (the identity matrix)
 and $cov(\M\y)-cov(\hat{\M}_{\beta}\y)$ is positive semi-definite
 for all matrix $\M$ subject to $\M\F=\I$.
 
 {\bf BLUP of $\eta$.} The BLUP of $\eta$ is the estimator $\hat{\eta}=\hat{\M}_{\eta}\y$, 
with $\hat{\M}_{\eta}$ (called the hat matrix of $\hat{\eta}$) being the matrix such that $\hat{\M}_{\eta}\F=0$
 and $cov(\M\y)-cov(\hat{\M}_{\eta}\y)$ is positive semi-definite
 for all matrix $\M$ subject to $\M\F=0$. 
 
We call, by convention, \emph{predictors} of a random variable
to distinguish them from \emph{estimators} of a deterministic parameter.
Henderson et al.(1959)\cite{Henderson3} showed that the BLUE and the BLUP are respectively 
 \ben
 \hat{\beta} &=& (\F^T(\R\sig_{\eta}\R^T+\sig_w)^{-1}\F)^{-1}\F^T(\R\sig_{\eta}\R^T+\sig_w)^{-1}\y\label{blue1959},\\
 \hat{\eta} &=& (\R^T\sig_N^{-1}\R+\sig_{\eta}^{-1})^{-1}[\R^T\sig_w^{-1}-\R^T\sig_{w}^{-1}\F(\F^T(\sig_{w}+\R\sig_{\eta}\R^T)^{-1}\F)^{-1}\nonumber \\
 &&\F^T(\sig_{w}+\R\sig_{\eta}\R^T)^{-1}]\y\label{blup1959}.
 \een 

Now we are able to give a Baysian interpretation of the hat matrix $\H_{\P_{pen}}(\lambda)$ (\ref{Plambdaformula}).  
Let $\P_1$ be an $n+3$ by $n+3-dim(N(\P_{pen}))$ matrix such that 
\ben\label{P1} 
\P_1^T\P_{pen}\P_1=\I_{(n+3-dim(N(\P_{pen})))\times (n+3-dim(N(\P_{pen})))},
\een 
respectively $\P_0$ an $n+3$ by $dim(N(\P_{pen}))$ matrix such that 
\ben\label{P0}  
\P_{pen}P_0=0,\quad  \mbox{and its columns form a basis of}\quad  N(\P_{pen}).
\een  
It follows that, for all vector $\lp{c} u_1\\ \p\\ u_{n+1}\rp$, there exist a unique $\beta\in \mathbb{R}^{dim(N(\P_{pen}))}$ and 
$\eta\in \mathbb{R}^{n+3-dim(N(\P_{pen}))}$ such that 
\ben 
\lp{c} u_1\\ \p\\ u_{n+1}\rp=\P_0\beta+\P_1\eta. 
\label{changeofvariable} 
\een 
Hence the model $\y=\p+\w$ becomes  
\be 
\y&=&\Pi P_0\beta+\Pi P_1\eta+\w\\
 &:=&\F\beta+\R\eta+\w.
 \ee 
We suppose that $\beta$ is the fixed effect and  $\eta$ is independent of the noise $\w$ and 
drawn from a centred distribution having the covariance matrix 
$\sigma_s^2\I_{n+3-dim(N(P))}$. 

Now, we are able to give our Bayesian interpretation. 
\begin{prop}  
The components $(\hat{\beta},\hat{\eta})$ of  
\be 
\arg\min\{\frac{\sigma_{w}^2}{\sigma_s^2}\|\eta\|^2+\|\y-(\F\beta+R\eta)\|^2:\quad 
\beta\in\mathbb{R}^{dim(N(\P_{pen}))}, \eta\in \mathbb{R}^{n+3-dim(N(\P_{pen}))}\}
\ee 
are respectively the BLUE of $\beta$ and the BLUP of $\eta$. 
Moreover, we have 
\be 
\F\hat{\beta}+\R\hat{\eta}&=&\Pi\H_{\P_{pen}}(\frac{\sigma_{w}^2}{\sigma_s^2})\lp{c}0\\ \y\\ 0\rp,\\
\P_0\hat{\beta}+\P_1\hat{\eta}&=&\H_{P_{pen}}(\frac{\sigma_{w}^2}{\sigma_s^2})\lp{c}0\\ \y\\ 0\rp.
\ee 
\end{prop} 

The proof is a consequence of the change of variable formula (\ref{changeofvariable})
and (\ref{blup1959}). 
See \cite{DRW} Proposition 2.2 for a similar proof.  

\begin{cor} Let $P=\C$ be the matrix (\ref{Cmatrix}) and 
$\P_0, \P_1$ be the corresponding matrices defined by (\ref{P0}), (\ref{P1}). We have 
\be 
\P_0\hat{\beta}&=&\lp{c}0\\ Lreg\y\\ 0\rp,\\
\P_1\hat{\eta}&=&\lp{c}0\\(\H(\frac{\sigma_{w}^2}{\sigma_s^2})-Lreg)\y\\ 0\rp.
\ee 
\end{cor}

{\bf Conclusion.} In this work we defined  a new basis of the set of $C^2$-cubic splines. We revisited the estimation of a natural cubic spline using Schoenberg-Reinsch result 
and we extended their result to the estimation of any $C^2$-cubic spline. 
We studied the choice of the smoothing parameter when the noise is deterministic or white throughout several criteria. We also gave a Bayesian interpretation of our estimators.

\section*{ Appendix 1: Construction of the matrices $\Q, \U, \V$} 
From (\ref{ctwo}) we can solve for $v_1, \ldots, v_{n}$ in terms 
of $u_1, \ldots, u_{n+1}$, i.e.  
\ben 
v_i=\frac{u_{i+1}-u_i}{h_i},\quad i=1, \ldots, n. 
\label{vu}
\een  
Now, using (\ref{cone}) and (\ref{czero}) we can solve $\q$ in terms of 
$\p$, $\u$. We get  
\be 
\q=\Q_1\p+\Q_2\u,
\ee 
where the $n$ by $n+1$ matrices  
\be 
\Q_1&=&[-\frac{1}{h_i}\e_i^T+\frac{1}{h_i}\e_{i+1}^T,\quad i=1,. . ., n],\\  
\Q_2&=&[-\frac{h_i}{3}\e_i^T-\frac{h_i}{6}\e_{i+1}^T, \quad i=1,. . ., n].   
\ee 
Here the column vectors $(\e_i: i=1, \ldots, n+1)$ denote the canonical basis of $\Rb^{n+1}$. 
If we plug $\q, \v=(v_1, \ldots,v_{n})^T$ in the continuity equations (\ref{czero}), then we get 
\ben 
\S\u=\Delta\p, 
\label{SDelta}
\een
where the $n-1$ by $n+1$ matrix 
\be 
\Delta=[\frac{e_i^T}{h_i}-(\frac{1}{h_i}+\frac{1}{h_{i+1}})e_{i+1}^T+\frac{1}{h_{i+1}}e_{i+2}^T,\quad i=1, \ldots, n-1].
\ee  
The $n-1$ by $n+1$ matrix
\ben 
\S=[\frac{h_i}{6}\e_i^T+\frac{h_i+h_{i+1}}{3}\e_{i+1}^T+\frac{h_{i+1}}{6}\e_{i+2}^T,\quad i=1,. . ., n-1]. 
\label{S} 
\een
Observe that the $n+3$ by $n+3$ matrix 
\be 
\S_{1,n+1}:=\lp{c} 
e_1^T\\
\S\\
e_{n+1}^T\rp
\ee 
is invertible. Hence, we can solve for $\u$ in terms 
of $(u_1,\p,u_{n+1})$ as follows: 
\ben 
\u=\U\lp{c} u_1\\ \p\\ u_{n+1}\rp
\label{u}
\een 
where the $n+1$ by $n+3$ matrix  
\ben 
\U&=&\S_{1,n+1}^{-1}\lp{c} 
e_1^T\\
0\quad \Delta\quad 0\\
e_{n+1}^T\rp.
\label{U}
\een 
The matrix $\U$ tells us, for all $i=1, \ldots, n+1$, that the second derivative $u_i$ is a linear combination of $u_1, u_{n+1}$ and 
$\p$ with the weight $\U=[u_{i,j}]$, i.e. 
\ben 
u_i=u_1u_{i,1}+\sum_{j=1}^{n+1}p_ju_{i,j+1}+u_{n+1}u_{i,n+3}.
\label{ui}
\een
The coefficient $u(i,1)$ is the weight of $u_1$, $(u_{i,j+1}, j=1, \ldots, n+1)$ are the weight of the
observations $(p_j, j=1, \ldots, n+1)$ respectively. The equality (\ref{ui}) tells us that the initial and the terminal second derivatives do not depend on 
the observations $\p$. We can show numerically and we will prove it rigorously (Proposition \ref{splinewithu1un}) that the mean weight 
of the obervations on each second derivative is equal to zero, i.e.,   
\be 
\sum_{j=1}^{n+1}u_{i,j+1}=0,\quad\forall\,i=1, \ldots, n+1.
\ee  

Now, we come back to the first and the third derivatives $\q, \v$.
We can solve for $\q$ and $\v$ in terms 
of $(u_1,\p,u_{n+1})$, i.e. 
\ben 
\q=\Q\lp{c} u_1\\ \p\\ u_{n+1}\rp,
\label{q}
\een 
with 
\be 
\Q=\lp{c} 
0\quad \Q_1\quad 0\rp+\Q_2\U.
 \ee
Similarly,  
\ben
\v=\V\lp{c} u_1\\ \p\\ u_{n+1}\rp,
\label{v}
\een 
with the $n$ by $n+3$ matrix 
\be 
\V=\tilde{V}\U,
\ee 
where the $n$ by $n+1$ matrix 
\be 
\tilde{V}=[-\frac{1}{h_i}e_i^T+\frac{1}{h_i}e_{i+1}^T,\quad i=1, \ldots, n].
\ee  
The matrices $\Q, \U, \V$ satisfy the equality 
\ben
\Q+diag(\frac{h}{2})\U+diag(\frac{h^2}{6})\V=[0\quad \D\quad 0],
\label{DE}
\een 
where the $n$ by $n+1$ matrix
\be 
\D=[\frac{1}{h_i}(e_{i+1}^T-e_i^T),\quad i=1, \ldots, n]. 
\ee 
\section*{Appendix 2: Analysis and interpretation of the matrices $\Q, \U, \V$} 
How to interpret the columns ?\\
$\bullet$ The derivatives up to order $3$, for $j=0, \ldots, n+2$,  of $\varphi_j$ at the knots $t_1, \ldots, t_{n+1}$ are respectively the $(j+1)$th columns of $\P, \Q, \U, \V$.  
Here $\P$ is the $n+1$ by $n+3$ matrix 
\be
\P=[\varphi_j(t_i): i=1,\ldots, n+1,\, j=0, \ldots, n+2].
\ee

How to interpret the rows ?\\
$\bullet$ The $i$th row of the matrices $\P, \Q, \U, \V$ represents respectively the row 
$(\varphi_j(t_i):\quad j=0, \ldots, n+2)$, $(\varphi_j'(t_i):\quad j=0, \ldots, n+2)$, $(\varphi_j''(t_i):\quad j=0, \ldots, n+2)$
and $(\varphi_j'''(t_i+):\quad j=0, \ldots, n+2)$.
 
How to interpret the basis elements?\\ 
$\bullet$ The natural cubic spline interpolating $(t_i,0),\,i=1, \ldots, n+1$ is the null map denoted by $s_0$. 
The natural cubic spline  interpolating $(t_j, 1), (t_i,0),\,i\neq j$ is equal to 
$\varphi_{j}$. It is the unique $C^2$-cubic spline such that, for $i=1, \ldots, n+1$,  
\be 
s(t_i)=\delta_i^j,\quad s'(t_i)=q_{i,j+1},\\
s''(t_i)=u_{i,j+1}, \quad s'''(t_i+)=v_{i,j+1}.
\ee   
It can be seen as the perturbation of the null map having the most important value at the knot $t_j$ (see Figure 1). 
A perturbation of $1$ at the knot $t_j$, i.e. $p_j\to p_j+1$, produces for each $t\in [t_1, t_{n+1}]\setminus\{t_j\}$
a perturbation $\varphi_j(t)\in (-1,1)$.\\ 
$\bullet$ The $C^2$-cubic spline $s$ interpolating $(t_i,0), i=1, \ldots, n+1$ and such that $s''(t_1)=1:=u_{1,1}$  is equal to 
$\varphi_{0}$.   
It is the unique $C^2$-cubic spline such that, for $i=1, \ldots, n+1$,  
\be 
s(t_i)=0,\quad s'(t_i)=q_{i,1},\\
s''(t_i)=u_{i,1}, \quad s'''(t_i+)=v_{i,1}.
\ee   
It can be seen as the perturbation of the null map having the most important second derivative at the knot $t_1$ (see Figure 3). 
A perturbation of $1$ on the second derivative at the knot $t_1$ produces for each 
$t\in (t_1, t_{n+1}]$
a perturbation $\varphi_0''(t)\in (-1,1)$.\\ 
$\bullet$ The $C^2$-cubic spline $s$ interpolating $(t_i,0), i=1, \ldots, n+1$ and such that $s''(t_{n+1})=1:=u_{1,n+1}$  is equal to 
$\varphi_{n+2}$.   
It is the unique $C^2$-cubic spline such that, for $i=1, \ldots, n+1$,  
\be 
s(t_i)=0,\quad s'(t_i)=q_{i,n+2},\\
s''(t_i)=u_{i,n+2}, \quad s'''(t_i+)=v_{i,n+2}.
\ee   
It can be seen as the perturbation of the null map having the most important second derivative at the knot $t_{n+1}$ (see Figure 3). 
A perturbation of $1$ on the second derivative at the knot $t_1$ produces for each $t\in [t_1, t_{n+1})$
a perturbation $\varphi_{n+2}''(t)\in (-1,1)$. 

\section*{Appendix 3: Calculus of the matrices $\C_{ij}$, $i,j=0, \ldots, 2$} 

We recall, for $i=1, \ldots, n$, $l=0, \ldots, n+2$ and $t\in (t_i,t_{i+1})$, that 
\be 
\varphi_l(t)=\delta_l(i)+q_{i,l+1}(t-t_i)+\frac{u_{i,l+1}}{2}(t-t_i)^2+\frac{v_{i,l+1}}{6}(t-t_i)^3,\\
\varphi_l'(t)=q_{i,l+1}+u_{i,l+1}(t-t_i)+\frac{v_{i,l+1}}{2}(t-t_i)^2,\\
\varphi_l''(t)=u_{i,l+1}+v_{i,l+1}(t-t_i).
\ee 

{\bf Calculus of $\C_{00}$.} We have, for $l,k=0, \ldots, n+2$, that 
\be 
\int_{t_1}^{t_{n+1}}\varphi_l(t)\varphi_k(t) dt=\sum_{i=1}^nh_i
\int_0^1\varphi_l(t,i)\varphi_k(t,i)dt, 
\ee 
where for $i=1, \ldots, n$, 
\be 
\varphi_l(t,i)&=&\delta_l(i)+q_{i,l+1}t+\frac{u_{i,l+1}}{2}t^2+\frac{v_{i,l+1}}{6}t^3\\
&:=&\sum_{p=0}^3a_p(l,i)t^p.
\ee 
Hence 
\be 
\varphi_l(t,i)\varphi_k(t,i)=\sum_{p=0}^6a_p(l,k,i)t^p,
\ee 
where 
\be 
a_p(l,k,i)=\sum_{p_1+p_2=p: p_1,p_2=0,1,2,3} a_{p_1}(l,i)a_{p_2}(k,i).
\ee 
It follows that 
\be 
\int_0^1\varphi_l(t,i)\varphi_k(t,i)dt=\sum_{p=0}^6\frac{a_p(l,k,i)}{p+1}, 
\ee 
and the matrix 
\be 
\C_{00}=[\sum_{i=1}^n\sum_{p=0}^6h_i\frac{a_p(l-1,k-1,i)}{p+1}:\quad l,k=1, \ldots, n+3]. 
\ee 

{\bf Calculus of $\C_{01}$.} We have, for $l,k=0, \ldots, n+2$, that 
\be 
\int_{t_1}^{t_{n+1}}\varphi_l(t)\varphi_k'(t) dt=\sum_{i=1}^nh_i
\int_0^1\varphi_l(t,i)\varphi_k'(t,i)dt, 
\ee 
where, for $i=1, \ldots, n$,  
\be 
\varphi_k'(t,i)&=&\sum_{p=1}^3 pa_p(k,i)t^{p-1}\\
&:=&\sum_{p=1}^3a_p^{(1)}(k,i)t^{p-1}.   
\ee 
Hence 
\be 
\varphi_l(t,i)\varphi_k'(t,i)=\sum_{p=1}^6 a_p^{(01)}(l,k,i)t^{p-1},
\ee 
where 
\be 
a_p^{(01)}(l,k,i)=\sum_{p_1+p_2=p: p_1=0,1,2,3, p_2=1,2,3} a_{p_1}(l,i)a_{p_2}^{(1)}(k,i).
\ee 
It follows that 
\be 
\int_0^1\varphi_l(t,i)\varphi_k'(t,i)dt=\sum_{p=1}^6\frac{a_p^{(01)}(l,k,i)}{p}, 
\ee 
and the matrix 
\be 
\C_{01}=[\sum_{i=1}^n\sum_{p=1}^6h_i\frac{a_p^{(01)}(l-1,k-1,i)}{p}: \quad l,k=1, \ldots, n+3].
\ee 

{\bf Calculus of $C_{02}$.} We have, for $l,k=0, \ldots, n+2$, that 
\be 
\int_{t_1}^{t_{n+1}}\varphi_l(t)\varphi_k''(t) dt=\sum_{i=1}^nh_i
\int_0^1\varphi_l(t,i)\varphi_k''(t,i)dt, 
\ee 
where, for $i=1, \ldots, n$,  
\be 
\varphi_k''(t,i)&=&\sum_{p=2}^3 p(p-1)a_p(k,i)t^{p-2}\\
&:=&\sum_{p=2}^3a_p^{(2)}(k,i)t^{p-2}.   
\ee 
Hence 
\be 
\varphi_l(t,i)\varphi_k''(t,i):=\sum_{p=2}^6a_p^{(2)}(l,k,i)t^{p-2},
\ee 
where 
\be 
a_p^{(02)}(l,k,i)=\sum_{p_1+p_2=p, p_1=0,1,2,3, p_2=2,3} a_{p_1}(l,i)a_{p_2}^{(2)}(k,i).
\ee 
It follows that 
\be 
\int_0^1\varphi_l(t,i)\varphi_k''(t,i)dt=\sum_{p=2}^6\frac{a_p^{(02)}(l,k,i)}{p-1}, 
\ee 
and
\be 
\C_{02}=[\sum_{i=1}^n\sum_{p=2}^6h_i\frac{a_p^{(02)}(l-1,k-1,i)}{p-1}:\quad l,k=1, \ldots, n+3].  
\ee 

{\bf Calculus of $\C_{11}$.} We have, for $l,k=0, \ldots, n+2$, that 
 \be 
 \varphi_l'(t,i)\varphi_k'(t,i)&=&\sum_{p_1,p_2=1}^3a_{p_1}^{(1)}(l,i)a_{p_2}^{(1)}(k,i)t^{p_1+p_2-2}\\
 &=&\sum_{p=2}^6a_p^{(11)}(l,k,i)t^{p-2},
 \ee 
 where 
 \be 
 a_p^{(11)}(l,k,i)=\sum_{p_1+p_2=p, p_1,p_2=1,2,3}a_{p_1}^{(1)}(l,i)a_{p_2}^{(1)}(k,i).
 \ee   
Hence, 
\be 
\int_{t_1}^{t_{n+1}}\varphi_l'(t)\varphi_k'(t) dt\\
&=&\sum_{i=1}^nh_i\int_0^1\varphi_l'(t,i)\varphi_k'(t,i)dt\\
&=&\sum_{i=1}^nh_i\sum_{p=2}^6\frac{a_p^{(11)}(l,k,i)}{p-1}, 
\ee 
and 
\be 
\C_{11}=[\sum_{i=1}^nh_i\sum_{p=2}^6\frac{a_p^{(11)}(l-1,k-1,i)}{p-1}:\quad l,k=1, \ldots, n+3].  
\ee 

{\bf Calculus of $\C_{12}$.} We have, for $l,k=0, \ldots, n+2$, that 

\be 
\varphi_l'(t,i)\varphi_k''(t,i)&=& 
\sum_{p_1=1}^3\sum_{p_2=2}^3a_{p_1}^{(1)}(l,i)t^{p_1-1}a_{p_2}^{(2)}(k,i)t^{p_2-2}\\
&=&\sum_{p=3}^6a_p^{(12)}(l,k,i)t^{p-3},       
\ee 
where 
\be 
a_{p}^{(12)}(l,k,i)=\sum_{p_1+p_2=p, p_1=1,2,3,p_2=2,3}a_{p_1}^{(1)}(l,i)a_{p_2}^{(2)}(k,i).
\ee 
Hence, 
\be 
\int_{t_1}^{t_{n+1}}\varphi_l'(t)\varphi_k''(t) dt\\
&=&\sum_{i=1}^nh_i\int_0^1\varphi_l'(t,i)\varphi_k''(t,i)dt\\
&=&\sum_{i=1}^n\sum_{p=3}^6h_i\frac{a_p^{(12)}(l,k,i)}{p-2}, 
\ee 
and 
\be 
\C_{12}=[\sum_{i=1}^nh_i\sum_{p=3}^6\frac{a_p^{(12)}(l-1,k-1,i)}{p-2}:\quad l,k=1, \ldots, n+3].
\ee 

\end{document}